\documentstyle[12pt]{article}
\input amssym.def
\input amssym
\begin{document}
\def \Z{\Bbb Z}
\def \C{\Bbb C}
\def \R{\Bbb R}
\def \Q{\Bbb Q}
\def \N{\Bbb N}
\def \wt{{\rm wt}}
\def \tr{{\rm tr}}
\def \span{{\rm span}}
\def \Res{{\rm Res}}
\def \Res{{\rm QRes}}
\def \End{{\rm End}}
\def \E{{\rm End}}
\def \Ind {{\rm Ind}}
\def \Irr {{\rm Irr}}
\def \Aut{{\rm Aut}}
\def \Hom{{\rm Hom}}
\def \mod{{\rm mod}}
\def \ann{{\rm Ann}}
\def \<{\langle}
\def \>{\rangle}
\def \t{\tau }
\def \a{\alpha }
\def \e{\epsilon }
\def \l{\lambda }
\def \L{\Lambda }
\def \g{\gamma}
\def \b{\beta }
\def \om{\omega }
\def \o{\omega }
\def \c{\chi}
\def \ch{\chi}
\def \cg{\chi_g}
\def \ag{\alpha_g}
\def \ah{\alpha_h}
\def \ph{\psi_h}
\def \be{\begin{equation}\label}
\def \ee{\end{equation}}
\def \bl{\begin{lem}\label}
\def \el{\end{lem}}
\def \bt{\begin{thm}\label}
\def \et{\end{thm}}
\def \bp{\begin{prop}\label}
\def \ep{\end{prop}}
\def \br{\begin{rem}\label}
\def \er{\end{rem}}
\def \bc{\begin{coro}\label}
\def \ec{\end{coro}}
\def \bd{\begin{de}\label}
\def \ed{\end{de}}
\def \pf{{\bf Proof. }}
\def \voa{{vertex operator algebra}}

\newtheorem{thm}{Theorem}[section]
\newtheorem{prop}[thm]{Proposition}
\newtheorem{coro}[thm]{Corollary}
\newtheorem{conj}[thm]{Conjecture}
\newtheorem{lem}[thm]{Lemma}
\newtheorem{rem}[thm]{Remark}
\newtheorem{de}[thm]{Definition}
\newtheorem{hy}[thm]{Hypothesis}
\makeatletter
\@addtoreset{equation}{section}
\def\theequation{\thesection.\arabic{equation}}
\makeatother
\makeatletter

\begin{center}{\Large \bf Regularity of rational vertex operator algebras}

\vspace{0.5cm}
Chongying Dong\footnote{Supported by NSF grant DMS-9303374 and a
research grant from the Committee on Research, UC Santa Cruz.},
Haisheng Li and Geoffrey Mason\footnote{Supported by NSF grant
DMS-9401272 and a research grant from the Committee on Research, UC
Santa Cruz.}\\ Department of Mathematics, University of California,
Santa Cruz, CA 95064
\end{center}

\section{Introduction}

Rational vertex operator algebras, which play a fundamental role
in rational conformal field theory (see [BPZ] and [MS]),
single out an important class of vertex
operator algebras. Most vertex operator algebras which have been studied
so far are rational vertex operator algebras. Familiar
examples include
the moonshine module  $V^{\natural}$ ([B], [FLM], [D2]), the vertex operator
algebras $V_L$ associated with positive definite even lattices $L$ ([B],
[FLM], [D1]),
the vertex operator algebras $L(l,0)$ associated with integrable
representations of affine Lie algebras [FZ] and the vertex operator algebras
$L(c_{p,q},0)$ associated with irreducible highest weight representations
for the discrete series of the  Virasoro algebra ([DMZ] and [W]).

A rational vertex operator algebra as studied in this paper is a
vertex operator algebra such that any {\em admissible} module is a
direct sum of simple ordinary modules (see Section 2). It is
natural to ask if such complete reducibility holds for an arbitrary
weak module (defined in Section 2). A rational vertex operator algebra
with this property is called a {\em regular} vertex operator algebra. One
motivation for studying such vertex operator algebras arises in trying to
understand the appearance of negative fusion rules (which are computed
by the Verlinde formula) for vertex operator algebras $L(l,0)$ for
certain rational $l$ (cf. [KS] and [MW]).

In this paper we give several sufficient conditions under which a rational
vertex operator algebra is regular. We prove that the
rational vertex operator algebras $V^{\natural},$ $L(l,0)$ for positive
integers $l,$ $L(c_{p,q},0)$ and $V_L$ for positive definite even lattices
$L$ are regular. Our result for $L(l,0)$
implies that any restricted integrable
module of level $l$ for the corresponding affine Lie
algebra is a direct sum of irreducible highest weight integrable modules.
This result is expected to be useful in comparing
the construction of tensor
product of modules for $L(l,0)$ in [F] based on Kazhdan-Lusztig's
approach [KL] with the construction of tensor product of modules [HL]
in this special case. We should remark that $V_L$ in general
is a vertex algebra in the sense of [DL] if $L$ is not positive definite.
In this case we establish the complete reducibility of any weak module.

Since the definition of vertex operator algebra is by now well-known, we
do not define vertex operator algebra in this paper. We refer the reader
to [FLM] and [FHL] for their elementary properties.
The reader can find the details of the constructions of $V^{\natural}$ and
$V_L$ in [FLM], and $L(l,0)$ and $L(c_{p,q},0)$ in [DMZ], [DL], [FLM], [FZ],
[L1] and [W].

The paper is organized as follows: In Section 2, after defining the
notion of weak module for a vertex operator algebra and the definition
of rational vertex operator algebra, we discuss the rational vertex
operator algebras $V^{\natural},$ $V_L,$ $L(l,0)$ and $L(c_{p,q},0).$
Section 3 is devoted to regular vertex operator algebras. We begin
this section with the definition of regular vertex operator
algebra. We show that the tensor product of regular vertex operator
algebras is also regular and that a rational vertex operator algebra
is regular under either of the assumptions (i) it contains a regular vertex
operator subalgebra, or (ii) any weak module contains a simple
ordinary module.
These
results are then used to prove that $V^{\natural},$ $V_L$ ($L$ is positive
definite), $L(l,0)$ and
$L(c_{p,q},0)$ are regular.  We also discuss the complete reducibility
of weak $V_L$-modules for an arbitrary even lattice $L.$ Based
on these results, we
conjecture that {\em any} rational vertex operator algebra is regular.

We thank Yi-Zhi Huang for pointing out a mistake in a prior version of
this paper.

\section{Rational vertex operator algebras}

Let $(V,Y,{\bf 1},\omega)$ be a vertex operator
algebra (cf. [B], [FHL] and [FLM]). A {\em weak module}
 $M$ for $V$ is a vector space equipped with a linear map
$$\begin{array}{l}
V\to (\mbox{End}\,M)[[z^{-1},z]]\label{map}\\
v\mapsto\displaystyle{ Y_M(v,z)=\sum_{n\in\Z}v_nz^{-n-1}\ \ \ (v_n\in
\mbox{End}\,M)}
\end{array}$$
(where for any vector space $W,$ we define $W[[z^{-1},z]]$ to be the vector
space of $W$-valued formal series in $z$)
satisfying the following conditions for $u,v\in V$,
$w\in M$:
\begin{eqnarray}
& &Y_M(v,z)=\sum_{n\in \Z}v_nz^{-n-1}\ \ \ \ \mbox{for}\ \
v\in V;\label{1/2}\\
& &v_nw=0\ \ \
\mbox{for}\ \ \ n\in \Z \ \ \mbox{sufficiently\ large};\label{vlw0}\\
& &Y_M({\bf 1},z)=1;\label{vacuum}
\end{eqnarray}
 \begin{equation}\label{jacobi}
\begin{array}{c}
\displaystyle{z^{-1}_0\delta\left(\frac{z_1-z_2}{z_0}\right)
Y_M(u,z_1)Y_M(v,z_2)-z^{-1}_0\delta\left(\frac{z_2-z_1}{-z_0}\right)
Y_M(v,z_2)Y_M(u,z_1)}\\
\displaystyle{=z_2^{-1}\delta\left(\frac{z_1-z_0}{z_2}\right)
Y_M(Y(u,z_0)v,z_2)}.
\end{array}
\end{equation}
\begin{equation}\label{vir}
[L(m),L(n)]=(m-n)L(m+n)+\frac{1}{12}(m^3-m)\delta_{m+n,0}(\mbox{rank}\,V)
\end{equation}
for $m, n\in {\Z},$ where
\begin{eqnarray}
& &L(n)=\omega_{n+1}\ \ \ \mbox{for}\ \ \ n\in{\Z}, \ \ \
\mbox{i.e.},\ \ \ Y_M(\omega,z)=\sum_{n\in{\Z}}L(n)z^{-n-2};\nonumber\\
& &\frac{d}{dz}Y_M(v,z)=Y_M(L(-1)v,z).\label{6.72}
\end{eqnarray}
This completes the definition. We denote this module by
$(M,Y_M)$ (or briefly by $M$).

\begin{de}\label{r2.2} An (ordinary) $V$-module is a  weak $V$-module which
carries a $\C$-grading
$$M=\coprod_{\lambda \in{\C}}M_{\lambda} $$
such that $\dim M_{\l}$ is finite and $M_{\l+n}=0$
for fixed $\l$ and $n\in {\Z}$ small enough. Moreover one requires that
$M_{\l}$ is the $\l$-eigenspace for $L(0):$
$$L(0)w=\l w=(\mbox{wt}\,w)w, \ \ \ w\in M_{\l}.$$
\end{de}

This definition is weaker than that of [FLM], for example, where the grading
on $M$ is taken to be rational. The extra flexibility attained by allowing
$\C$-gradings is important $-$ see for example [DLM1] and [Z].

We observe some redundancy in the definition of weak module:
\begin{lem}\label{r2.1} Relations (\ref{vir}) and (\ref{6.72}) in the
definition of weak module are consequences of (\ref{1/2})-(\ref{jacobi}).
\end{lem}

\pf To establish (\ref{6.72}) note that
$L(-1)u=L(-1)u_{-1}{\bf 1}=u_{-2}{\bf 1}$ for $u\in V.$ Then
\begin{eqnarray}
& &\ \ \ Y_{M}(L(-1)u,z_{2})\nonumber\\
& &=Y_{M}(u_{-2}{\bf 1},z_{2})\nonumber\\
& &={\rm Res}_{z_{0}}z_{0}^{-2}Y_{M}(Y(u,z_{0}){\bf 1},z_{2})\nonumber\\
& &={\rm Res}_{z_{0}}{\rm Res}_{z_{1}}z_{0}^{-2}
\left(z^{-1}_0\delta\left(\frac{z_1-z_2}{z_0}\right)
Y_M(u,z_1)Y_M({\bf 1},z_2)\right.\nonumber\\
& &\ \ \ \ \ \ \left.-z^{-1}_0\delta\left(\frac{z_2-z_1}{-z_0}\right)
Y_M({\bf 1},z_2)Y_M(u,z_1)\right)\nonumber\\
& &={\rm Res}_{z_{0}}{\rm
Res}_{z_{1}}z_{0}^{-2}z_2^{-1}\delta\left(\frac{z_1-z_0}{z_2}\right)
Y_{M}(u,z_{1})\nonumber\\
& &={\rm Res}_{z_{0}}{\rm
Res}_{z_{1}}z_{0}^{-2}z_1^{-1}\delta\left(\frac{z_2+z_0}{z_1}\right)
Y_{M}(u,z_{2}+z_{0})\nonumber\\
& &={\rm Res}_{z_{0}}z_{0}^{-2}Y_{M}(u,z_{2}+z_{0})\nonumber\\
& &={\rm Res}_{z_{0}}z_{0}^{-2}e^{z_0\frac{d}{d z_2}}Y_{M}(u,z_{2})\nonumber\\
& &={d\over dz_{2}}Y_{M}(u,z_{2}).
\end{eqnarray}
This establishes (\ref{6.72}), and together with
(\ref{jacobi}) and  $Y(\omega,z_0)\omega=\frac{1}{2}({\rm rank}V)z_0^{-4}
+2\omega z_0^{-2}+L(-1)\omega z_0^{-1}+$ regular terms we can easily deduce
(\ref{vir}). \ \ \ $\Box$

Thus we may just use  (\ref{1/2})-(\ref{jacobi}) as the axioms
for a weak $V$-module.

\begin{de}\label{d2.2} An {\em admissible} $V$-module is
a  weak $V$-module $M$ which carries a
${\Z}_{+}$-grading
$$M=\coprod_{n\in {\Z}_{+}}M(n)$$
($\Z_+$ is the set all nonnegative integers) satisfying
the following condition: if $r, m\in {\Z} ,n\in {\Z}_{+}$ and $a\in V_{r}$
then
\begin{eqnarray}
a_{m}M(n)\subseteq M(r+n-m-1).\label{2.7}
\end{eqnarray}

We call an admissible $V$-module $M$ {\em simple}
in case $0$ and $M$ are the only $\Z_+$-graded submodules.

$V$ is called {\em rational} if every admissible
$V$-module is a direct sum of simple admissible $V$-modules. That is,
we have complete reducibility of admissible $V$-modules.
\end{de}

\begin{rem}\label{r2.4} (i) Note that any ordinary $V$-module is admissible.

(ii) It is proved in [DLM1] that
if $V$ is rational then conversely, every simple admissible $V$-module is an
ordinary module. Moreover $V$ has only a finite number of inequivalent
simple modules.

(iii) Zhu's definition of rational vertex operator algebra $V$ is as follows
[Z]: (a) all admissible $V$-module are completely reducible, (b) each simple
admissible $V$-module is an ordinary $V$-module, (c) $V$
only has finitely many inequivalent simple modules. Thanks to (ii), Zhu's
definition of rational thus coincides with our own.
\end{rem}

We next introduce a certain category
$\cal{O}$ of admissible $V$-modules in analogy with the well-known
category $\cal O$ of Bernstein-Gelfand-Gelfand. First some
notation: for any weak $V$-module $M$ we set for $h\in \C:$
$$M_{h}=\{m\in M|
(L(0)-h)^km=0\ {\rm for\ some} \ k\in\Z_{+}\}.$$
So $M_h$ is a {\em generalized} eigenspace for $L(0),$ and  in particular
$M_h$  is  the $h$-eigenspace for $L(0)$ if $L(0)$ is a semisimple operator.

Now define $\cal{O}$ to be the category of weak $V$-modules $M$ satisfying
the following two conditions:

(1) $L(0)$ is locally finite in the sense that if $m\in M$ then
there is a finite-dimensional $L(0)$-stable subspace of $M$
which contains $m.$

(2) There are $h_{1}, \cdots, h_{k}\in\C$ such that $$M=\oplus
_{i=1}^{k}\oplus _{n\in {\Z}_{+}}M_{n+h_{i}}.$$ These are the {\em
objects} of $\cal O.$ Morphisms may be taken to be $V$-module
homomorphisms, though we will not make use of them in the sequel.

\begin{rem}\label{r2.5} (i) Any weak $V$-module which belongs to $\cal O$
is necessarily admissible: a $\Z_+$-grading obtains by defining
$M(n)=\oplus_{i=1}^kM_{h_i+n}.$ Condition (\ref{2.7}) follows in the
usual way.

(ii) Suppose that $M$ is a weak $V$-module and that $W$ is a weak
$V$-submodule of $M.$ Then $M$ lies in $\cal O$ if, and only if,
both $W$ and $M/W$ lie in $\cal O.$

(iii) If $V$ is rational, any weak $V$-module
in $\cal{O}$ is a direct sum of simple $V$-modules (use Remark \ref{r2.4}
(ii)).
\end{rem}

Next we briefly discuss some familiar examples of rational vertex
operator algebras. The reader is referred to the references for notation
and the details of the constructions.

(1) Let $L$ be an even lattice and $V_{L}$ the corresponding
vertex algebra (see [B], [DL] and [FLM]). It is proved in [D1] that
if $L$ is positive definite then
$V_{L}$ is rational and its simple modules are parametrized by
$L'/L$ where $L'$ is the dual lattice of $L.$

(2) Let ${\frak g}$ be a finite-dimensional simple Lie algebra with a Cartan
subalgebra ${\frak h}$ and $\hat \frak g=\C[t,t^{-1}]\otimes \frak g\oplus\C c$
the corresponding affine Lie algebra. Fix a positive integer $l.$
 Then any $\lambda\in \frak h^*$ can be viewed as a linear
form on $\C c\oplus \frak h\subset \hat \frak g$ by sending $c$ to $l.$
Let us denote the
corresponding irreducible highest weight module for $\hat\frak g$ by
$L(l,\lambda).$ Then $L(\ell,0)$ is a rational vertex operator algebra
([DL], [FZ], [L1]).

(3) Let $c$ and $h$ be two complex numbers and let $L(c,h)$ be the lowest
weight irreducible
module for the Virasoro algebra with central charge $c$ and lowest weight $h$.
Then $L(c,0)$ has
a natural vertex operator algebra structure (cf. [FZ]).
Moreover, $L(\ell,0)$ is rational if, and only if,
$c=c_{p,q}=1-\frac{6(p-q)^{2}}{pq}$ for $p,q\in \{2,3,4,\cdots\}$ and
$p, q$ are relatively prime (see [DMZ] and [W]).

(4) Let $V^{\natural}$ be the moonshine module vertex operator algebra
constructed by Frenkel, Meurman and Lepowsky [FLM] (see also [B]).
It is established in [D2] that $V^{\natural}$ is {\em holomorphic} in the
sense that $V^{\natural}$ is rational and the only simple
module is $V^{\natural}$ itself.

(5) Let $V^{1},\cdots, V^{k}$ be vertex operator algebras.
Then $V=\otimes_{i=1}^{k}V^{i}$ is
a vertex operator algebra of rank $\sum_{i=1}^{k}{\rm rank} V^{i}$ and
any simple $V$-module $M$ is isomorphic to a tensor product module
$M^{1}\otimes \cdots \otimes M^{k}$ for some simple $V^{i}$-module $M^i$
[FHL].
Furthermore, $V^{1}\otimes V^{2}\otimes \cdots\otimes V^{k}$ is rational
if, and only if, each $V^i$ is rational [DMZ].

(6) Let $V^{1}, \cdots, V^{k}$ be vertex operator algebras of the same
rank. Then $\oplus_{i=1}^{k}V^{k}$ is a vertex operator algebra [FHL].
It is clear that
$\oplus_{i=1}^{k}V^{i}$ is rational if each $V^{i}$ is rational. The vacuum
space of the resulting vertex operator algebra is not one-dimensional,  and
we will not consider this particular example further in this
paper.

\section{Regular vertex operator algebras}

In Section 2 we made use of the complete reducibility of admissible modules in
order to
define rational vertex operator algebras.
The study of complete reducibility of an arbitrary weak module
for a rational vertex operator algebra
leads us to the following notion of regular vertex operator algebra:

\begin{de}\label{d2.4}
A vertex operator algebra $V$ is said to be {\em regular} if any weak
$V$-module $M$ is a direct
sum of simple ordinary $V$-modules.
\end{de}

\begin{rem}\label{r3.2}  A regular vertex operator algebra $V$ is necessarily
rational. Indeed if $M$ is a weak $V$-module then, being a direct sum
of ordinary simple $V$-module, it is admissible (Remark \ref{r2.4} (i)) and,
for the same reason, a direct sum of simple admissibles.
\end{rem}

The main result of the present paper is to show that the rational
vertex operator algebras in examples (1)-(4) of Section 2 are each
regular. First we have some general results which will be useful
later.

\begin{prop}\label{p2.5}
Let $V^{1}, \cdots, V^{k}$ be regular vertex operator algebras. Then
$V=V^{1}\otimes V^{2}\otimes \cdots \otimes V^{k}$ is regular.
\end{prop}

{\bf Proof.} Let $M$ be a weak $V$-module. For each $1\le i\le k$, we may
regard $V^{i}$ as a vertex operator subalgebra with a different Virasoro
element. Then
$M$ is a weak $V^{i}$-module. Since $V^i$ is regular, $M$ is a direct
sum of simple ordinary $V^{i}$-modules. Note that
there are only finitely many simple $V^{i}$-modules up to
equivalence. Thus $M$ is a $V^{i}$-module in category $\cal{O}$ of weak
$V^i$-modules.
Denote the generators of the Virasoro algebra of $V^i$ by $L_i(n).$
Then $L(0)=L_{1}(0)+\cdots +L_{k}(0)$ and $L_{i}(0)$'s commute with each
other. This implies that $M$ is in category $\cal{O}$ of weak
$V$-modules. So $M$ is completely reducible by Remarks \ref{r3.2} and
\ref{r2.5} (iii). $\;\;\;\;\Box$

\begin{prop}\label{p2.6}
Let $V$ be a rational vertex operator algebra such that there is a regular
vertex operator subalgebra $U$
with the same Virasoro element $\omega$. Then $V$ is regular.
\end{prop}

{\bf Proof.} Let $M$ be a weak $V$-module. Then $M$ is a weak $U$-module,
 so that $M$ is a direct
 sum of simple ordinary $U$-modules. Thus $L(0)$ acts semisimply on $M$.
Let $W^{1},\cdots, W^{k}$ be all simple $U$-modules up to equivalence.
Then we can write $M=\sum_{i=1}^{k}\oplus_{n\in {\Z}_{+}}M_{n+h_{i}}$,
where $h_{i}$ is the lowest weight of $W^{i}$.
Thus $M$ is in the category $\cal{O}$ for $V.$
The complete reducibility of $M$ follows immediately as
$V$ is rational.$\;\;\;\;\Box$

\begin{prop}\label{p2.8}
Let $V$ be a rational vertex operator algebra such that any nonzero
weak $V$-module contains a simple ordinary $V$-submodule. Then $V$ is
regular.
\end{prop}

{\bf Proof.} Let $M$ be any weak $V$-module and let $W$ be the sum of all
simple ordinary submodules.
We have to prove that $W=M$. If $M\ne W$, $M/W$ is a nonzero weak
$V$-module so that by assumption, there is a simple ordinary
$V$-submodule $M^{1}/W$ of $M/W$.
Then both $W$ and $M^{1}/W$ are in the category $\cal{O}$.
Thus $M^{1}$ is in the category $\cal{O}$ and $M^{1}$ is a direct sum of
simple ordinary $V$-modules. This contradicts
the choice of $W$. $\;\;\;\;\Box$

Now we are ready to
show that the
vertex operator algebras $L(l,0),$ $L(c_{p,q},0),$ $V^{\natural}$ and
$V_L$  are regular.

Recall from example (2) that $\frak g$ is a finite-dimensional simple Lie
algebra with Cartan subalgebra $\frak h;$ $L(l,0)$ is a vertex operator
algebra. We shall denote the corresponding root system by $\Delta.$

\begin{lem}\label{l2.9} There is a basis $\{a^{1},\cdots, a^{m}\}$ for ${\frak
g}$ such that for $1\leq i,j\leq m$ we have
\begin{eqnarray}\label{eb}
[Y(a^{i},z_{1}),Y(a^{i},z_{2})]=0\;\;\mbox{ and }Y(a^{i},z)^{3\ell+1}=0
\end{eqnarray}
as operators on $L(\ell,0).$
\end{lem}

{\bf Proof.} Let $\alpha\in\Delta$ and $e\in {\frak g}_{\alpha}$.
If ${\frak g}$ is of
type $A, D,$ or $ E$, it is proved in [MP1] and [DL] that $Y(e,z)^{\ell+1}=0$.
In general, it is proved in [L1] and [MP2] that $Y(e,z)^{3\ell+1}=0$.
It is well known (cf. [H], [K1]) that there are elements
 $e_{\alpha}\in {\frak g}_{\alpha}, f_{\alpha}\in {\frak g}_{-\alpha},
h_{\alpha}\in \frak h$
which linearly span a subalgebra naturally isomorphic to $sl_{2}$. Set
$\sigma_{\alpha}=e^{(e_{\alpha})_{0}}$
where $(e_{\alpha})_0$ is the component operator of $Y(e_{\alpha},z)$
(cf. (\ref{1/2})) corresponding to  $z^{-1}.$
Then $\sigma_{\alpha}$ is an automorphism of the vertex
operator algebra $L(\ell,0)$ (see Chap. 11 of [FLM]).  A straightforward
calculation gives
$\sigma_{\alpha}(f_{\alpha})=f_{\alpha}+h_{\alpha}-2e_{\alpha}$. Since
$e_{\alpha}, f_{\alpha}, h_{\alpha}$ form a basis of ${\frak g}$ for
$\alpha\in \Delta,$
$e_{\alpha}, f_{\alpha}, \sigma_{\alpha}(f_{\alpha})$ also form a basis of
${\frak g}$. It is clear that this basis satisfies condition (\ref{eb}).
$\;\;\;\;\Box$

\bt{t2.10}
Let $\ell$ be a positive integer. Then
the vertex operator algebra $L(\ell,0)$ is regular.
\et

{\bf Proof.} By Proposition \ref{p2.8}, it is enough to prove that any nonzero
weak
$L(\ell,0)$-module $M$ contains a simple $L(\ell,0)$-module. This
will be established in three steps.

{\bf Claim 1:} {\em There exists a nonzero $u\in M$ such that
$(t{\C}[t]\otimes {\frak g})u=0$.} Set $\frak g(n)=t^n\otimes \frak g$ for
$n\ne 0.$
For any nonzero $u\in M$, by the definition of a weak module,
${\frak g}(n)u=0$ for sufficiently large
$n.$ So $(t{\C}[t]\otimes {\frak g})u$ is finite-dimensional.
For any $u\in M$, we define $d(u)=\dim (t{\C}[t]\otimes {\frak g})u$.
If there is a $0\ne u\in M$ such that $d(u)=0$, then
$(t{\C}[t]\otimes {\frak g})u=0$. Suppose that $d(u)>0$ for any $0\ne u\in M$.
Take $0\ne u\in M$ such that $d(u)$ is minimal.

Let $a^{i}$ ($1\le i\le m$) be a basis of ${\frak g}$ satisfying
condition (\ref{eb}) and let $k$ be the positive integer such that
${\frak g}(k)u\ne 0$ and ${\frak g}(n)u=0$ whenever $n>k$. By
definition of $k$, $a^{i}(k)u\ne 0$ for some $1\le i\le m$. Since
$Y(a^{i},z)^{3\ell+1}=0$, by Proposition 13.16 in [DL],
$Y_{M}(a^{i},z)^{3\ell+1}=0$. Extracting the coefficient of
$z^{-(k+1)(3\ell+1)}$ from $Y_{M}(a^{i},z)^{3\ell+1}u=0$ we obtain
$(a^{i}_{k})^{3\ell+1}u=0$.

Let $r$ be a nonnegative integer  such that $(a^{i}_{k})^{r}u\ne 0$
and $(a^{i}_{k})^{r+1}u=0$. Set $v=(a^{i}_{k})^{r}u$. We will obtain
a contradiction by showing that
$d(v)<d(u)$.
First we prove that if $a_{n}u=0$ for some $a\in {\frak g}, 1\le n\in {\Z}$,
then $a_{n}v=0$. In the following we will show by induction
on $m$ that $a_n(a_k^i)^mu=0$ for
any $a\in \frak g$ and $m\in \Z$ nonnegative. If $m=0$ this is immediate by the
choice of $u.$ Now assume that the result holds for $m.$
Since $[a,a^{i}]_{k+n}u=0$ (from the definition of $k$) and $a_{n}u=0,$ by the
induction assumption that $a_n(a_k^i)^mu=0$
we have:
\begin{eqnarray}
[a,a^{i}]_{k+n}(a^{i}_{k})^{m}u=0,\;a_{n}(a^{i}_{k})^{m}u=0.
\end{eqnarray}
Thus
\begin{eqnarray}
&
&a_{n}(a^{i}_{k})^{m+1}u=[a_{n},a^{i}_{k}](a^{i}_{k})^{m}u+a^{i}_{k}a_{n}(a^{i}_{k})^{m}u\nonumber\\
& &\ \ \ \
=[a,a^{i}]_{k+n}(a^{i}_{k})^{m}u+a^{i}_{k}a_{n}(a^{i}_{k})^{m}u\nonumber\\
& &\ \ \ \ =0,
\end{eqnarray}
as required. In particular, we see that
$a_{n}v=a_{n}(a^{i}_{k})^{r}u=0$.
Therefore, $d(v)\le d(u)$. Since $a^{i}_{k}v=0$ and $a^{i}_{k}u\ne 0$, we have
$d(v)< d(u)$.

{\bf Claim 2:} {\em There is a nonzero  $u\in M$ such that ${\frak g}(n)u=0$
for $n>0$ and ${\frak g}_+u=0$ where $\frak g_+=\sum_{\alpha\in \Delta_+}\frak
g_{\alpha}$
for a fixed positive root system $\Delta_+.$}

 Set
\begin{eqnarray}
\Omega(M)=\{u\in M| (t{\C}[t]\otimes {\frak g})u=0\}.
\end{eqnarray}
Then $\Omega(M)$ is a nonzero ${\frak g}$-submodule of $M$ by Claim 1.
Let $0\ne e_{\theta}\in {\frak g}_{\theta}$
where $\theta$ is the longest positive root in $\Delta$. Then
$Y_{M}(e_{\theta},z)^{\ell+1}=0$ (see [DL] and [FZ]). Extracting the
coefficient of $z^{-\ell-1}$ from
$Y_{M}(e_{\theta},z)^{\ell+1}\Omega(M)=0$,
we obtain $e_{\theta}^{\ell+1}\Omega(M)=0$.
By Proposition 5.1.2 of [L1] $\Omega(M)$ is a direct sum of finite-dimensional
irreducible
${\frak g}$-modules. Then any highest weight vector for ${\frak g}$ in
$\Omega(M)$ meets our need.

{\bf Claim 3:} {\em Any lowest weight vector for $\hat{{\frak g}}$ in $M$
generates a simple
$L(\ell,0)$-module.} Let $u$ be a lowest weight vector for $\hat{{\frak g}}$ in
$M$.
Extracting the constant term from
$Y_{M}(e_{\theta},z)^{\ell+1}u=0$, we obtain $(e_{\theta})_{-1}^{\ell+1}u=0$.
Then $u$ generates an integrable highest weight $\hat{{\frak g}}$-module. It
follows from [K1] that
$u$ generates an irreducible $\hat{{\frak g}}$-module of level $\ell$.
Since any submodule of $M$
for the affine Lie algebra is a submodule of $M$ for $L(\ell,0)$,
such $u$ generates a simple
$L(\ell,0)$-module.$\;\;\;\;\Box$

\begin{rem}\label{r2.3'} This theorem has been proved
in [DLM2] under the assumption that $t\C[t]\otimes \frak h$ acts locally
nilpotently on any weak module. See Proposition 5.6 in [DLM2].
\end{rem}

\begin{rem}
Recall that a $\hat{{\frak g}}$-module $M$ is called  {\em restricted }
(cf. [K1]) if for any $u\in M$, there is an integer $k$ such that
 $(t^{n}\otimes {\frak g})u=0$ for $n>k;$
$M$ is called an {\em integrable} module
if the Chevalley generators $e_i,f_i$ of $\hat{{\frak g}}$
act locally finitely on $M$ [K2] (note that in the definition of
integrable module, we do not assume that the action of $\frak h$
is semisimple). At affine Lie algebra level, Proposition \ref{p2.8}
and Theorem \ref{t2.10}
essentially assert
that any restricted integrable $\hat{{\frak g}}$-module is a direct sum of
irreducible highest weight integrable $\hat{{\frak g}}$-modules.
\end{rem}

Next we turn our attention to the vertex operator algebras $L(c_{p,q},0).$

First we recall some results from [DL]. Let $V$ be a vertex operator
algebra and $M$ be a weak $V$-module. Then for any $u,v\in V$ we have
\begin{eqnarray}\label{ef}
Y(u_{-1}v,z)=Y(u,z)^{-}Y(v,z)+Y(v,z)Y(u,z)^{+},
\end{eqnarray}
where
\begin{eqnarray}
Y(u,z)^{+}=\sum_{n\ge 0}u_{n}z^{-n-1}, \; Y(u,z)^{-}=\sum_{n<0}u_{n}z^{-n-1}.
\end{eqnarray}
Note that $Y(u,z)=Y(u,z)^{+}+Y(u,z)^{-}$ and that $Y(u,z)^{+}$ (reps.
$Y(u,z)^{-}$) involves only
nonpositive (resp. nonnegative) powers of $z$.

For convenience we will write $c=c_{p,q}.$
For any nonnegative integer $n$, we set
$\omega^{(n)}={1\over n!}L(-1)^{n}\omega$. Then $L(-n-2)=\omega^{(n)}_{-1}.$
We need the following lemmas.

\begin{lem}\label{l2.12}
Let $M$ be a weak $L(c,0)$-module and $u\in M.$ Let
$k$ be a positive integer such that
$L(k)u\ne 0$ and that $L(n)u=0$ whenever $n>k$. Then for any nonnegative
integers $n_{1},..., n_{r}$ the lowest power of $z$ in
$Y_{M}(\omega^{(n_{1})}_{-1}\cdots \omega^{(n_{r})}_{-1}{\bf 1},z)u$
(in the sense that the coefficients of $z^m$ is zero
whenever $m$ is smaller than the lowest power) is
$-r(k+2)-n_{1}-\cdots -n_{r}$ with  coefficient
$\displaystyle{\prod_{i=1}^{r}{-k-2\choose n_{i}}L(k)^{r}u}.$
\end{lem}

{\bf Proof.} We prove this lemma by induction on $r$. If $r=1$, we have:
\begin{eqnarray}
& &Y_M(\omega^{(n_{1})}_{-1}{\bf 1},z)u=Y_M(\omega^{(n_{1})},z)u\nonumber\\
& &\ \ \ \ \ =\frac{1}{n_1!}\left({d\over
dz}\right)^{n_{1}}Y_M(\omega,z)u\nonumber\\
& &\ \ \ \ \ =\sum_{n\in {\Z}}{-n-2\choose n_{1}}z^{-n-2-n_{1}}L(n)u
\nonumber\\
& &\ \ \ \ \ =\sum_{n\le k}{-n-2\choose n_{1}}z^{-n-2-n_{1}}L(n)u.
\end{eqnarray}
Then the lowest power of $z$ is $-(k+2)-n_{1}$ with a coefficient
${-k-2\choose n_{1}}L(k)u$. That is, the
lemma holds for $r=1$.

Suppose that this lemma holds for some positive integer $r$.
By formula (\ref{ef}) we have:
\begin{eqnarray}
& &Y_{M}(\omega^{(n_{1})}_{-1}\cdots
\omega^{(n_{r})}_{-1}\omega^{(n_{r+1})}_{-1}{\bf
1},z)u=Y_{M}(\omega^{(n_{1})},z)^{-}Y_{M}(\omega^{(n_{2})}_{-1}\cdots
\omega^{(n_{r+1})}_{-1}{\bf 1},z)
u\nonumber\\
& &\ \ \ \ \ \ +Y_{M}(\omega^{(n_{2})}_{-1}\cdots \omega^{(n_{r+1})}_{-1}{\bf
1},z)
Y_{M}(\omega^{(n_{1})},z)^{+}u.\label{a2.15}
\end{eqnarray}
Since $Y_{M}(\omega^{(n_{1})},z)^{-}$ involves only nonnegative powers of $z$,
it follows from
the inductive assumption, the lowest power of $z$ in the first term of the
right hand side of (\ref{a2.15}) is
$-r(k+2)-n_{2}-\cdots -n_{r+1}$. It is easy to observe that for
any $v$ in the algebra,
$$Y_M(L(-1)v,z)^+=\frac{d}{dz}Y_M(v,z)^+.$$
Thus
\begin{eqnarray*}
& &\ \ \ \ Y_{M}(\omega^{(n_{2})}_{-1}\cdots \omega^{(n_{r+1})}_{-1}{\bf 1},z)
Y_{M}(\omega^{(n_{1})},z)^{+}u\\
& &=\frac{1}{n_1!}Y_{M}(\omega^{(n_{2})}_{-1}\cdots \omega^{(n_{r+1})}_{-1}
{\bf 1},z)\left(\frac{d}{dz}\right)^{n_1}Y_{M}(\omega,z)^{+}u\\
& &=\sum_{n=-1}^k{-n-2\choose n_1}Y_{M}(\omega^{(n_{2})}_{-1}\cdots
\omega^{(n_{r+1})}_{-1}{\bf 1},z)L(n)uz^{-n-2-n_1}\\
& &=\sum_{n=0}^k{-n-2\choose n_1}Y_{M}(\omega^{(n_{2})}_{-1}\cdots
\omega^{(n_{r+1})}_{-1}{\bf 1},z)L(n)uz^{-n-2-n_1}\\
& &\ \ \ \ \ +{-1\choose n_1}L(-1)Y_{M}(\omega^{(n_{2})}_{-1}\cdots
\omega^{(n_{r+1})}_{-1}{\bf 1},z)uz^{-1-n_1}\\
& &\ \ \ \ \ +{-1\choose
n_1}\left(\frac{d}{dz}Y_{M}(\omega^{(n_{2})}_{-1}\cdots
\omega^{(n_{r+1})}_{-1}{\bf 1},z)\right)uz^{-1-n_1}.
\end{eqnarray*}
Note that $L(m)L(n)u=(m-n)L(m+n)u+L(n)L(m)u=0$ for $0\leq n\leq k$ and
$m>k.$ Applying the inductive hypothesis to $L(n)u$ we see that the lowest
 power of $z$ in
$$\sum_{n=0}^k{-n-2\choose n_1}Y_{M}(\omega^{(n_{2})}_{-1}\cdots
\omega^{(n_{r+1})}_{-1}{\bf 1},z)L(n)uz^{-n-2-n_1}$$
is $-(r+1)(k+2)-n_1-\cdots -n_{r+1}$ with coefficient
$\displaystyle{\prod_{i=1}^{r+1}{-k-2\choose n_i}L(k)^{r+1}u.}$  Also by the
induction assumption
the  lowest  power of $z$ in  $L(-1)Y_{M}(\omega^{(n_{2})}_{-1}\cdots
\omega^{(n_{r+1})}_{-1}{\bf 1},z)uz^{-1-n_1}$ is
$-r(k+2)-n_1-\cdots -n_{r+1}-1$ and the  lowest  power of $z$ in
$\left(\frac{d}{dz}Y_{M}(\omega^{(n_{2})}_{-1}\cdots
\omega^{(n_{r+1})}_{-1}{\bf 1},z)\right)uz^{-1-n_1}$
is $-r(k+2)-n_1-\cdots -n_{r+1}-2.$
Thus  the lowest power of
$z$ in the second
term of right hand side of (\ref{a2.15}) is $-(r+1)(k+2)-n_{1}-\cdots -n_{r+1}$
with coefficient
$\displaystyle{\prod_{i=1}^{r+1}\left(\begin{array}{c}-k-2\\n_{i}\end{array}\right)
L(k)^{r+1}u},$ as desired. $\;\;\;\;\Box$

\bigskip

Let $V$ be a vertex operator algebra and let $A(V)$  (which is a certain
quotient space
of $V$ modulo a subspace $O(V)$) be the corresponding
associative algebra defined in [Z]. We refer the reader to [Z] for details.
Recall from [DLM1] or [L2] that for any weak $V$-module $M,$
$\Omega(M)$ consists of vectors $u\in M$ such that
$a_{m}u=0$ for any homogeneous element $a\in V$ and for any $m>{\rm wt}\,a-1$.
In other words,
$u\in \Omega(M)$ if and only if $z^{m}Y_{M}(a,z)u\in M[[z]]$ for any
homogeneous element $a\in V$
and for any $m>{\rm wt}a-1$. The following result can be found in [DLM1], [L2]
and [Z].

\begin{lem}\label{al} (1) $\omega+O(V)$ is in the center of $A(V).$

(2) $A(V)$ is semisimple if $V$ is rational.

(3) $\Omega(M)$ is an $A(V)$-module under the action  $a+O(V)\mapsto
a_{\wt\,a-1}$ for homogeneous $a\in V.$
\end{lem}

Now we take $V=L(c,0).$
Set $\bar{\Omega}(M)=\{ u\in M|L(n)u=0\;\mbox{ for any }n>0\}$. Then it is
clear that
$\Omega(M)\subseteq \bar{\Omega}(M)$.

\begin{lem}\label{l2.13}
Let $M$ be a weak $L(c,0)$-module. Then $\Omega(M)=\bar{\Omega}(M)$.
\end{lem}

{\bf Proof.} It suffices to prove that $a_{m}u=0$ for any $u\in
\bar{\Omega}(M)$ and for any
homogeneous element $a\in L(c,0)$ whenever $m>{\rm wt}a-1$. We shall prove
this by induction on the weight of $a$.
If ${\rm wt}a=0$, $a={\bf 1}$. Since ${\bf 1}_{m}=0$ for $m\ge 0$, there is
nothing to prove.
Suppose that $a_{m}u=0$ for any homogeneous element $a\in L(c,0)$ of weight
less than $n$ and for any
$m>{\rm wt}a-1$. Let $b\in L(c,0)$ be a homogeneous element of weight $n$ and
let $m\in {\Z}$ such that $m>{\rm wt}b-1$.

Let $a\in L(c,0)$ be any homogeneous element of weight less than $n$, let $k$
be any positive integer
and let $m>{\wt}\,(L(-k)a)-1$ $(={\rm wt}\,a+k-1)$. Then from the Jacobi
identity (\ref{jacobi}) we have:
\begin{eqnarray}
& &(L(-k)a)_{m}u={\rm Res}_{z_{0}}{\rm
Res}_{z_{2}}z_{0}^{1-k}z_{2}^{m}Y_{M}(Y(\omega,z_{0})a,z_{2})u\nonumber\\
& &\ \ ={\rm Res}_{z_{1}}{\rm Res}_{z_{0}}{\rm
Res}_{z_{2}}z_{0}^{1-k}z_{2}^{m}\cdot \nonumber\\
& &\ \ \ \cdot\left(
z_{0}^{-1}\delta\left(\frac{z_{1}-z_{2}}{z_{0}}\right)Y_{M}(\omega,z_{1})Y_{M}(a,z_{2})u
-z_{0}^{-1}\delta\left(\frac{z_{2}-z_{1}}{-z_{0}}\right)Y_{M}(a,z_{2})Y_{M}(\omega,z_{1})u\right)
\nonumber\\
& &\ \ ={\rm Res}_{z_{1}}{\rm
Res}_{z_{2}}(z_{1}-z_{2})^{1-k}z_{2}^{m}Y_{M}(\omega,z_{1})Y_{M}(a,z_{2})u
\nonumber\\
& &\ \ \ -{\rm Res}_{z_{1}}{\rm
Res}_{z_{2}}(-z_{2}+z_{1})^{1-k}z_{2}^{m}Y_{M}(a,z_{2})Y_{M}(\omega,z_{1})u.
\end{eqnarray}
Since $m>{\rm wt}\,(L(-k)a)-1={\rm wt}a+k-1>{\rm wt}a-1$, we have:
$${\rm Res}_{z_{1}}{\rm
Res}_{z_{2}}(z_{1}-z_{2})^{1-k}z_{2}^{m}Y_{M}(\omega,z_{1})Y_{M}(a,z_{2})u=0.$$
For the second term, we have:
\begin{eqnarray}
& &\ \ \ -{\rm Res}_{z_{1}}{\rm
Res}_{z_{2}}(-z_{2}+z_{1})^{1-k}z_{2}^{m}Y_{M}(a,z_{2})Y_{M}(\omega,z_{1})u
\nonumber\\
& &=-{\rm Res}_{z_{2}}(-1)^{1-k}z_{2}^{m+1-k}Y_{M}(a,z_{2})L(-1)u
-{\rm Res}_{z_{2}}(-1)^{-k}(1-k)z_{2}^{m-k}Y_{M}(a,z_{2})L(0)u\nonumber\\
& &=-{\rm Res}_{z_{2}}(-1)^{1-k}z_{2}^{m+1-k}L(-1)Y_{M}(a,z_{2})u
+{\rm Res}_{z_{2}}(-1)^{1-k}z_{2}^{m+1-k}{d\over
dz_{2}}Y_{M}(a,z_{2})u\nonumber\\
& &\ \ \ -{\rm
Res}_{z_{2}}(-1)^{-k}(1-k)z_{2}^{m-k}Y_{M}(a,z_{2})L(0)u\nonumber\\
& &=-{\rm Res}_{z_{2}}(-1)^{1-k}z_{2}^{m+1-k}L(-1)Y_{M}(a,z_{2})u
-{\rm Res}_{z_{2}}(-1)^{1-k}(m+1-k)z_{2}^{m-k}Y_{M}(a,z_{2})u\nonumber\\
& &\ \ \ -{\rm
Res}_{z_{2}}(-1)^{-k}(1-k)z_{2}^{m-k}Y_{M}(a,z_{2})L(0)u\nonumber\\
&
&=(-1)^kL(-1)a_{m+1-k}u+(-1)(m+1-k)a_{m-k}u+(-1)^{k-1}a_{m-k}L(0)u.\label{2.17}
\end{eqnarray}
Since $L(0)u\in \bar{\Omega}(M)$ and $m-k>{\rm wt}a-1$ all the three terms
in (\ref{2.17}) are zero by the inductive hypothesis.  Thus $(L(-k)a)_mu=0.$
Note that $b$ is a linear combination of all
$L(-k)a$, where ${\rm wt}a<n$ and $k$ is a positive integer.
This shows  $b_{m}\bar{\Omega}(M)=0$ for $m>{\rm wt}\,b-1,$ as desired.
$\;\;\;\;\Box$

Now we are in a position to prove
\bt{t2.11}
The vertex operator algebra $L(c,0)$ associated with the lowest weight
irreducible module for
the Virasoro algebra with central charge $c=c_{p,q}$ is regular.
\et

{\bf Proof:} By Proposition \ref{p2.8}, it is enough to prove that any
nonzero weak $L(c,0)$-module $M$ contains a simple $L(c,0)$-module.

{\bf Claim 1:} {\em The space $\Omega(M)$ is not zero.}
For any $0\ne u\in M$, we define $l(u)$ to be the integer $k$ such that
$L(k)u\ne 0$ and $L(n)u=0$
whenever $n>k$. Since $L(n)u\ne 0$ for some $n$ (because $c\ne 0$),
$l(u)$ is well-defined.
Suppose that $\Omega(M)=0$. Then by Lemma \ref{l2.13}
$l(u)\ge 1$ for any $0\ne u\in M$. Let $0\ne u\in M$ such that $l(u)=k$ is
minimal.
It is well known [FF] that there are two singular vectors in the
Verma module $M(c,0)$ for the Virasoro algebra. One singular
vector is $L(-1){\bf 1}$ and the other is:
\begin{eqnarray}
v=L(-2)^{pq}{\bf 1}+\sum a_{n_{1},\cdots, n_{r}}\omega^{(n_{1})}_{-1}\cdots
\omega^{(n_{r})}_{-1}{\bf 1},
\end{eqnarray}
where the sum is over some $(n_{1},\cdots, n_{r})\in {\Z}_{+}^{r}$
such that $2pq=2r+n_{1}+\cdots +n_{r}$ and $n_{1}+\cdots +n_{r}\ne
0$. By Lemma \ref{l2.12}, the lowest power of $z$ in
$$Y_{M}(L(-2)^{pq}{\bf 1},z)u=Y_{M}((\omega_{-1})^{m}{\bf 1},z)u$$ is
$-pq(k+2)$ with $L(k)^{pq}u$ as its coefficient and the lowest power
of $z$ in
$$Y_{M}(\omega^{(n_{1})}_{-1}\cdots \omega^{(n_{r})}_{-1}{\bf
1},z)u$$
 is greater than $-pq(k+2)$ for any nonnegative integers
$n_{1},\cdots, n_{r}$ such that $2pq=2r+n_{1}+\cdots +n_{r}$ and
$n_{1}+\cdots +n_{r}\ne 0$. Thus the coefficient of $z^{-pq(k+2)}$ in
$Y_{M}(v,z)u$ is $L(k)^{pq}u$. Since $v=0$ in $L(c,0)$ we have
$Y_{M}(v,z)=0$.  In particular the coefficient $L(k)^{pq}u$
of $z^{-pq(k+2)}$ in $Y_M(v,z)$ is zero. Let $s$ be the nonnegative
integer such that $L(k)^{s}u\ne 0$ and $L(k)^{s+1}u=0$ and set
$u'=L(k)^{s}u$. Then it is clear that $l(u')<l(u)$. This is a
contradiction.

{\bf Claim 2:} {\em Any weak $L(c,0)$-module $M$ contains a simple ordinary
$L(c,0)$-module.} Since $L(c,0)$ is rational, Lemma \ref{al} tells us that
$A(L(c,0))$ is semisimple and that
the central element $\omega+O(L(c,0))$ acts
semisimply on $\Omega(M)$ as $L(0).$  Since $\Omega(M)$ is nonzero by Claim
1 we can take $0\ne u\in \Omega(M)$  such that $L(0)u=hu$ where $h\in \C.$
Again since $L(c,0)$ is rational [W], $u$ generates a simple (ordinary)
$L(c,0)$-module. The proof is complete.
$\;\;\;\;\Box$

\begin{coro}\label{cm}
The moonshine module vertex operator algebra $V^{\natural}$ is regular.
\end{coro}

{\bf Proof.} From [DMZ], $V^{\natural}$ contains $L({1\over
2},0)^{\otimes 48}$ as a vertex operator subalgebra. Then the result follows
from Theorem \ref{t2.11}, Propositions \ref{p2.5} and \ref{p2.6}.
$\;\;\;\;\Box$

Finally we discuss the complete reducibility of
weak $V_L$-modules for an even lattice $L.$ We refer the reader to
[FLM] and [D1] for the construction of $V_L$ and related notations.

Let $M$ be any weak $V_L$-module. Define the vacuum space
$$\Omega_M=\{u\in M| \alpha(i)u=0\ {\rm for}\ \a\in L,\ i>0\}.$$
\begin{lem}\label{p2.14}
Let $L$ be an even lattice. Then for any
weak $V_{L}$-module $M,$ $\Omega_M\ne 0.$
\end{lem}

{\bf Proof.} For $u\in M$ then $A_u=\span\{\a(n)u|\alpha\in L, n>0\}$ is
finite-dimensional as $\a(n)u=0$ if $n$ is sufficiently large and as
the rank of $L$ is finite. Set $d(u)=\dim A_u.$ Note that $d(u)=0$ if
and only if $u\in \Omega_M.$ So it is enough to show that $d(u)=0$ for
some nonzero $u\in M.$ Assume this is false, and
take $0\ne u\in M$ such that $d(u)$ is
minimal.

Let $k$ be the smallest positive integer such that $\a(k)u\ne 0$ and
$\b(n)u=0$ whenever $n>k$ for some $\a\in L$ and all $\b\in L.$

Let $a\in \hat L$ such that $\bar a=\a.$ Then from the
formula (3.4) of [D1] we have
$$\frac{d}{dz}Y(\iota(a),z)=Y(L(-1)\iota(a),z)=Y(\alpha(-1)\iota(a),z)
=\alpha(z)^-Y(\iota(a),z)+Y(\iota(a),z)\alpha(z)^+$$
where
$$\alpha(z)^-=\sum_{n<0}\alpha(n)z^{-n-1},\ \
\alpha(z)^+=\sum_{n\geq 0}\alpha(n)z^{-n-1}.$$
Clearly the submodule generated by $u$ is not zero. Note that the vertex
algebra $V_L$ is simple (see [D1]).
By Proposition 11.9
of [DL], $Y(\iota(a),z)u\ne 0.$ Let $r$ be an integer such that
$\iota(a)_{r+m}u=0$ and $\iota(a)_{r}u\ne 0$ for any positive integer
$m.$ Thus the lowest power of $z$ in
$$\frac{d}{dz}Y(\iota(a),z)u=-\sum_{m\leq r}(m+1)\iota(a)_muz^{-m-2}$$
is at most $-r-2.$ It is obvious that the lowest power of $z$ in $
\alpha(z)^-Y(\iota(a),z)u$ is at most $-r-1.$

Use the following commutator formula which is a result from the
Jacobi identity
\begin{equation}\label{gr}
[\beta(m),\iota(a)_n]=\<\a,\b\>\iota(a)_{m+n}
\end{equation}
to obtain
$$\iota(a)_m\alpha(n)u=-\<\a,\a\>\iota(a)_{m+n}u+\alpha(n)\iota(a)_mu=0$$
if $m>r$ and $n\geq 0.$ This gives
$$Y(\iota(a),z)\alpha(z)^+u=
\sum_{m\leq r}\sum_{n=0}^k\iota(a)_m\a(n)z^{-m-n-2}.$$
Thus the coefficient $\iota(a)_{r}\a(k)u$ of $z^{-r-k-2}$ in the formula
above is zero as $k$ is positive. This shows by (\ref{gr}) again
that $\alpha(n)\iota(a)_ru=0$ for any positive integer greater than or equal
to $k.$

Note from (\ref{gr}) that if $\b(m)u=0$ for positive $m$ then
$\b(m)\iota(a)_ru=0.$ Thus $d(\iota(a)_ru)< d(u).$ This is a contradiction.
$\;\;\;\;\Box$

\bt{tvl}
Let $L$ be an even lattice. Then any weak $V_L$-module
is completely reducible and any simple weak $V_L$-module is isomorphic
to $V_{L+\beta}$ for some $\b$ in the dual lattice of $L.$
In particular, $V_L$ is regular if $L$ is positive definite.
\et

{\bf Proof.} By Lemma \ref{p2.14}, $\Omega_M\ne 0$ for a weak $V_L$-module
$M$.
It is proved in [D1] that if $M$ is also simple then it
is necessarily isomorphic to
$V_{L+\beta}$ for some $\b$ in the dual lattice of $L.$  So
it remains to show
the complete reducibility of any weak $V_L$-module $M.$

Let $W$ be the sum of all simple submodules of $M.$ Assume that
$M'=M/W$ is not zero. Then $\Omega_{M'}\ne 0.$ It is essentially proved
in [D1] that $M'$ contains a simple module $W^1/W$ (here $W^1$ is
a weak $V_L$-submodule of $M$ which contains $W$) generated
by $w^1+W$ where $w^1$ is a common eigenvector for the operators
$\a(0)$ for $\a\in L.$ It follows from the proof of Theorem 3.1 of [D1]
that the submodule of $M$ generated by $w^1$ is simple. Thus
$W^1$ is a sum of certain simple submodules of $V_L.$  This is
a contradiction.
$\;\;\;\;\Box$

At this point we have proved that almost all known rational vertex
operator algebras are regular. We conclude this paper by presenting
the following conjecture:

\begin{conj}
Any rational vertex operator algebra is regular.
\end{conj}

\end{document}